\documentclass[letterpaper, 10 pt, conference]{ieeeconf}  \IEEEoverridecommandlockouts                            
\let\proof\relax

\usepackage{bm}
\usepackage{amssymb}
\usepackage{makeidx}         
\usepackage{amsmath}
\usepackage{multicol}        
\usepackage[bottom]{footmisc}
\usepackage{graphicx}
\usepackage{color}
\usepackage{amsfonts}
\usepackage[latin1]{inputenc}
\usepackage{epstopdf}
\usepackage[usenames,dvipsnames]{xcolor}
\usepackage[ruled,vlined]{algorithm2e}
\usepackage{soul}
\usepackage{hyperref}

\usepackage{subcaption}	\SetKwInput{KwInput}{Input}
\SetKwInput{KwOutput}{Output}

\newcommand{\eq}{\begin{equation}}
\newcommand{\eeq}{\end{equation}}
\newcommand{\eqn}{\begin{eqnarray}}
\newcommand{\eeqn}{\end{eqnarray}}

\newcommand{\bsea}{\begin{subeqnarray}}
\newcommand{\esea}{\end{subeqnarray}}
\newcommand{\nn}{\nonumber}

\newcommand{\al}[1]{\begin{align}  #1 \end{align}}

\def\bmat{\left[ \begin{array}}
\def\emat{\end{array} \right]}

\newcommand{\Bc}{ \mathcal{B}}

\newcommand{\Dc}{ \mathcal{D}}

\newcommand{\Jc}{ \mathcal{J}}

\newcommand{\Rs}{ \mathbb{R}}

\def\qed{\hfill \vrule height 7pt width 7pt depth 0pt \smallskip}

\newcounter{pippo}

\newtheorem{remark}{Remark}[section]
\newtheorem{teor}{Theorem}[section]
\newtheorem{corr}{Corollary}[section]
\newtheorem{propo}{Proposition}[section]
\newtheorem{lemm}{Lemma}[section]
\newtheorem{exam}{Example}
\newtheorem{probl}[pippo]{Problem}
\newtheorem{defn}{Definition}[section]
\newcommand{\proof}{\noindent {\bf Proof. }}
\newcommand{\teo}{\begin{teor}}
\newcommand{\eteo}{\end{teor}}
\newcommand{\cor}{\begin{corr}}
\newcommand{\ecor}{\end{corr}}
\newcommand{\prop}{\begin{propo}}
\newcommand{\eprop}{\end{propo}}
\newcommand{\lem}{\begin{lemm}}
\newcommand{\elem}{\end{lemm}}
\newcommand{\ex}{\begin{exam}}
\newcommand{\eex}{\end{exam}}
\newcommand{\pb}{\begin{probl}}
\newcommand{\epb}{\end{probl}}
\newcommand{\df}{\begin{defn}}
\newcommand{\edf}{\end{defn}}
\newcommand{\aprop}{\begin{apropo}}
\newcommand{\eaprop}{\end{apropo}}
\newcommand{\alem}{\begin{alemm}}
\newcommand{\ealem}{\end{alemm}}
\newcommand{\rem}{\begin{remark}}
\newcommand{\erem}{\end{remark}}

\title{\textbf{Image compression by means of the multidimensional \\ circulant covariance extension problem -- Revisited} }

\author{Tommaso Benciolini, Tommaso Grigoletto and Mattia Zorzi \thanks{T. Benciolini, T. Grigoletto  M. Zorzi are with the Department of Information Engineering, University of Padova, Via Gradenigo 6/B, 35131 Padova, Italy. Emails: {\tt\small tommaso.benciolini@studenti.unipd.it}, {\tt\small tommaso.grigoletto@studenti.unipd.it}, {\tt\small zorzimat@dei.unipd.it}}}

\begin{document}

\maketitle
\thispagestyle{empty}
\pagestyle{empty}

\begin{abstract} We revisit the image compression problem using the framework introduced by Ringh, Karlsson and Lindquist. More precisely, we explore the possibility to consider a family of objective functions and a different way to design the prior in the corresponding multidimensional circulant covariance extension problem. The latter leads to refined compression paradigms.  \end{abstract}

\section{Introduction}
Image compression is a fundamental task to reduce the cost for storage or transmission. In the present paper we consider the image compression framework proposed in \cite{7403052, SIAM_IM_COMPR} which has been considered also in further extensions, see \cite{ringh2017further,ringh2018multidimensional,zhu2016identification}. Here, an image is characterized by a positive function defined over a 2-dimensional grid. The compressed image is constituted by a   finite set of moments of such a function. Then, the image reconstruction is performed by solving a moment matching problem: find a 2-dimensional positive function matching the given moments and maximizing a suitable objective function. The latter can be understood as a multidimensional circulant covariance extension problem. Indeed, such a positive function can be understood as the power spectral density of a periodic random Markov field and its moments correspond to the covariance lags. From this problem it is also possible to derive spectral estimation methods which can be extended also to the case of nonperiodic Markov random fields, \cite{georgiou2006relative,karlsson2016multidimensional,8660416,zhu2019fusion}. Moreover, in the case that the domain of spectral density boils down to a 1-dimensional grid, we obtain a circulant covariance extension problem corresponding to periodic stationary stochastic processes \cite{carli2011covariance,lindquist2013multivariate,lindquist2013circulant,ringh2015fast}. 

The covariance extension problem has been formerly studied to design high resolution spectral estimators for stationary stochastic processes, \cite{byrnes2000new,pavon2013geometry,8706531,zorzi2012estimation} including the approximate moments matching case \cite{enqvist2007approximative,ZORZI_ECC20}. Within this framework, an important aspect is that it is possible to take as objective function a pseudo-distance (or divergence) between the spectral density to be estimated and a given spectral density, called prior. The latter represents the a priori information that we have about the process, \cite{Hellinger_Ferrante_Pavon,FERRANTE_TIME_AND_SPECTRAL_2012,georgiou2003kullback}. In plain words, it is possible to select as optimal solution the closest one to the prior and matching the moments. A second important aspect is that it is possible to take a divergence family as objective function, \cite{ALPHA,BETAPRED,BETA,DUAL,zorzi2019graphical}. The latter leads to a family of solutions characterized by an integer parameter. 

The aim of this this paper is to revisit the image compression approach in \cite{7403052} by considering 
the aforementioned aspects. More precisely, in the corresponding multidimensional circulant covariance extension problem, we use the Alpha divergence family as objective function as well as new ideas to construct the prior in order to refine the compression paradigm.
 
The outline of the paper is as follows: in Section \ref{sect_back} we review the image compression framework proposed in \cite{7403052}; Sections \ref{sect_obj} and  \ref{sect_prior}  are devoted to the choice of the objective function and the prior, respectively; in Section \ref{sect_concl} we draw the conclusions.
 
\section{Brackground} \label{sect_back} 
A grayscale image can be represented by a nonnegative matrix whose entries, representing the pixels, take value in the interval $[0,1]$. Consider an image with $p_1\times p_2$ pixels. Then, the corresponding matrix is $X\in\mathbb R^{p_1\times p_2}$.  An equivalent representation is given by the function $\Phi(\zeta_{\bm \ell})$ with $\zeta_{\bm \ell}=[\, e^{i\ell_1\frac{2\pi}{N_1}}\; e^{i\ell_2\frac{2\pi}{N_2}} \, ]^T$, $\bm \ell=[\,\ell_1\;\ell_2\,]^T$, $N=[N_1\, N_2]^T\in\mathbb Z^2$ and  $\bm \ell\in\mathbb Z_N^2:= \{\, (\ell_1,\ell_2) \hbox{ s.t. } 0\leq \ell_1\leq N_1-1,\; 0\leq \ell_2\leq N_2-1\,\}$. More precisely, the relation between $X$ and $\Phi$ is as follows. We extend the image by symmetric mirroring: 
\al{Y_{i,j} = \begin{cases}
 X_{i,j} &\text{ if } i\leq p_1,j\leq p_2\\
 X_{i,(p_2-j+1)} &\text{ if } i\leq p_1,j> p_2\\
 X_{(p_1-i+1),j} &\text{ if } i> p_1,j\leq p_2\\
 X_{(p_1-i+1),(p_2-j+1)} &\text{ if } i> p_1,j> p_2\\
 \end{cases}}
with $Y\in\mathbb R^{N_1\times N_2}$, $X_{i,j}$ denotes the entry in position $(i,j)$ of $X$, $N_1=2(p_1-1)$ and $N_2=2(p_2-1)$. Then, we have $\Phi(\zeta_{\bm \ell})=\exp(Y_{\ell_1+1,\ell_2+1})$. Notice that $\Phi(\zeta_{\bm \ell})>0$, that is $\Phi$ is a multidimensional positive function.

Let $\Phi_o(\zeta_{\bm \ell})$ represent the original image. The latter  can be compressed by means of moments:
\begin{equation}\label{moment_cond}
c_{\bm k}=\frac{1}{|N|}\sum_{\bm \ell\in\mathbb Z_N^2} \zeta_{\bm \ell}^{\bm k} \Phi_o(\zeta_{\bm \ell}), \; \; \bm k\in\Lambda
\end{equation}
where $\zeta_{\bm \ell}^{\bm k}:=\zeta_{\ell_1}^{k_1}\zeta_{\ell_2}^{k_2}$, $\bm k=[\,k_1\; k_2\,]^T$, $\Lambda:=\{ \, [\, k_1\; k_2\,]^T\in \mathbb Z^2 \hbox{ s.t. } |k_1|\leq n_1,\; |k_2|\leq n_2\,\}$ with $n_1\ll N_1$, $n_2\ll N_2$ and $|N|=N_1N_2$. Accordingly, the moments set $\{c_{\bm k}, \, \bm k\in\Lambda \}$ represents the compressed image. Since $c_{\bm k} \in\mathbb R$ and $c_{\bm k}=c_{-\bm k}$, then the moments set is characterized by $(n_1+1)(n_2+1)$ parameters. Notice that $n_1$ and $n_2$ are fixed by the user and characterize the compression rate. The previous compression strategy is effective if we are able to extract a ``good'' approximation of $\Phi_o$ from $\{c_{\bm k}, \, \bm k\in\Lambda \}$. The latter can be understood as a multidimensional covariance extension problem: $\Phi$ represents the power spectral density of a periodic random Markov field and $c_{\bm k}$'s are the corresponding covariance lags. Clearly, given $\{c_{\bm k}, \, \bm k\in\Lambda\}$, there are infinite positive multidimensional functions $\Phi$ satisfying the moments constraint  in (\ref{moment_cond}). In order to choose one of the aforementioned solutions, we have to solve the following optimization problem:

\begin{equation}\label{opt}
\hat \Phi=\underset{\Phi>0}{\mathrm{argmin}}  \; \mathcal J(\Phi) \quad \hbox{ s.t. (\ref{moment_cond}) holds} 
\end{equation}
where $\mathcal J$ is a suitable cost function which guarantees that: 1) the optimization problem does admit solution and such a solution, say $\hat \Phi$, is unique; 2) $\hat \Phi$ is a good approximation of the original image $\Phi_o$. A natural choice is the entropy functional (with changed sign) $\Jc(\Phi)=-|N|^{-1}\sum_{\bm \ell \in\mathbb Z_N^2} \log \Phi(\zeta_{\bm \ell})$. Indeed, the latter choses the most ``flat'' solution matching the moments constraint. In \cite{7403052}, the authors proposed a refined method which considers the generalized entropy functional (with changed sign):
\al{\label{def_entropy}\Jc(\Phi)=-\frac{1}{|N|}\sum_{\bm \ell \in\mathbb Z_N^2} \Psi(\zeta_{\bm \ell})\log \Phi(\zeta_{\bm \ell})}
where $\Psi$ is a multidimensional trigonometric polynomial associated with the index set $\Lambda$. The latter is characterized in such a way that $\hat \Phi$ satisfies both the moment constraint in (\ref{moment_cond}) as well as the cepstral moments with index in $\Lambda\setminus\{0\}$. However, the existence of such a solution is not guaranteed. Such an issue is addressed by considering a regularized version of the problem whose solution approximately fulfils the cepstral matching, see the former work \cite{PER_CEPSTRAL}. It is worth noting that $\Psi$, hereafter called prior, embeds some information about the image that has to be reconstructed. In the following sections we explore the possibility to consider various objective functions in (\ref{opt}) as well as different ways to design the prior $\Psi$.

\section{Choice of the objective function}\label{sect_obj}
In this section we assume that the prior $\Psi(\zeta_{\bm \ell} )>0$, embedding some information about the image, is fixed. Then, a natural choice for the objective function in problem (\ref{opt}) is $\Jc(\Phi)=\Dc(\Phi\|\Psi)$ where $\Dc(\Phi\|\Psi)$ is a divergence between two positive functions $\Phi(\zeta_{\bm \ell}),\Psi(\zeta_{\bm \ell})>0$ such that $\Dc(\Phi\|\Psi)\geq 0$ and equality hold if and only if $\Phi=\Psi$. 

We choose as $\Dc$ the Aplha divergence family with the parametrization $\alpha=1-\frac{1}{\nu}$ and $\nu\in\mathbb N$, \cite{ALPHA}: for $1<\nu<\infty$, we have
\al{&\Dc_{\nu}(\Phi\|\Psi)=\nonumber \\ &\frac{1}{|N|}\sum_{\bm \ell \in\mathbb Z_N^2} \frac{\nu^2}{1-\nu}\Phi(\zeta_{\bm \ell})^{\frac{\nu-1}{\nu}}\Psi(\zeta_{\bm \ell})^{\frac{1}{\nu}}+\nu \Phi(\zeta_{\bm \ell})+\frac{\nu}{\nu-1}\Psi(\zeta_{\bm \ell});\nn} for $\nu=1$ and $\nu=\infty$, we have, respectively,
{\small \al{\Dc_{1}(\Phi\|\Psi)&=\frac{1}{|N|}\sum_{\bm \ell \in\mathbb Z_N^2}\Psi(\zeta_{\bm \ell})\log\left(\frac{\Psi(\zeta_{\bm \ell})}{\Phi(\zeta_{\bm \ell})}\right)- \Psi(\zeta_{\bm \ell})+\Phi(\zeta_{\bm \ell}),\nonumber\\
\Dc_{\infty}(\Phi\|\Psi)&=\Dc_{1}(\Psi\|\Phi)\nn.}}%
Therefore, we aim to solve the following problem

\begin{equation}\label{opt2}
\hat \Phi_\nu=\underset{\Phi>0}{\mathrm{argmin}}  \; \Dc_\nu (\Phi\|\Psi) \quad
 \hbox{ s.t. (\ref{moment_cond}) holds} 
\end{equation}
where $\Psi(\zeta_{\bm \ell})>0$ has been fixed and $1\leq \nu\leq \infty$. It is worth noting that the reconstructed image $\hat \Phi_\nu$ depends on the parameter $\nu$. 

In what follows we aim to show that problem (\ref{opt2}) does admit a unique solution through the duality theory. The latter also provides an efficient algorithm to solve the problem. We start with the case $1<\nu<\infty$. First, notice that (\ref{opt2}) is equivalent to 
\begin{equation}\label{opt3}
\hat \Phi_\nu=\underset{\Phi>0}{\mathrm{argmin}}  \; \frac{1}{\nu}\Dc_\nu (\Phi\|\Psi) \quad
\hbox{ s.t. (\ref{moment_cond}) holds.} 
\end{equation}

The Lagrangian is 
\al{L&(\Phi,Q)=\frac{1}{\nu}\Dc_\nu(\Phi\|\Psi)+\sum_{ \bm k\in\Lambda} q_{\bm k}\left(\frac{1}{|N|}\sum_{\bm \ell\in\mathbb Z_N^2} \zeta_{\bm \ell}^{\bm k} \Phi(\zeta_{\bm \ell})-c_{\bm k}\right)\nn\\
&=\frac{1}{|N|}\sum_{\bm \ell\in\mathbb Z_N^2}  \frac{\nu}{1-\nu} \Phi (\zeta_{\bm \ell})^{\frac{\nu-1}{\nu}}\Psi(\zeta_{\bm \ell})^{\frac{1}{\nu}}+Q(\zeta_{\bm \ell})\Phi (\zeta_{\bm \ell})\nn\\ & \hspace{0.5cm}-\sum_{\bm k\in\Lambda} q_{\bm k} c_{\bm k}+b}%
where $Q(\zeta_{\bm \ell})=\sum_{ \bm k\in\Lambda} q_{\bm k}\zeta^{\bm k}_{\bm \ell}$ is the Lagrange multiplier and $b$ is a constant term not depending on $\Phi$ and $Q$. Notice that $b$ contains also the constant term
\al{\frac{1}{|N|}\sum_{\bm \ell\in\mathbb Z_N^2} \Phi(\zeta_{\bm \ell}) = c_{\bm 0}\label{eqn:c0}}
fixed by the moment constraints in (\ref{moment_cond}). It is not difficult to see that $L(\cdot, Q)$ is strictly convex for $\Phi(\zeta_{\bm \ell})>0$. The first variation of $L(\cdot,Q)$ along the direction $\delta \Phi(\zeta_{\bm \ell})$ is
\al{\delta L(\Phi,Q;\delta\Phi)=\frac{1}{|N|}\sum_{\bm \ell\in\mathbb Z_N^2}(Q(\zeta_{\bm \ell})- \Phi(\zeta_{\bm \ell})^{-\frac{1}{\nu}}\Psi(\zeta_{\bm \ell})^{\frac{1}{\nu}})\delta \Phi(\zeta_{\bm \ell}). \nn}
Accordingly, the minimum of $L(\cdot ,Q)$ must satisfy the stationarity condition $\delta L(\Phi,Q;\delta \Phi)=0$ for any $\delta \Phi(\zeta_{\bm \ell})$. The latter implies$Q(\zeta_{\bm \ell})-\Phi(\zeta_{\bm \ell})^{-\frac{1}{\nu}}\Psi(\zeta_{\bm \ell})^{\frac{1}{\nu}}=0$. 
 The point of minimum is thus
$\Phi(\zeta_{\bm \ell})={\Psi (\zeta_{\bm \ell})  }/{Q(\zeta_{\bm \ell})^{\nu}}$ with $Q(\zeta_{\bm \ell})>0$, in order to satisfy $\Phi(\zeta_{\bm \ell})>0$. Observe that this is the same result of~\cite{ALPHA} up to the fact that we considered the term in~\eqref{eqn:c0} as constant, thus not depending on the optimization variable. Substituting the point of minimum in the Lagrangian, we obtain the dual functional (with changed sign)
\al{J_\nu(Q)= \frac{1}{|N|}\sum_{\bm \ell\in\mathbb Z_N^2} \frac{1}{\nu-1} \frac{\Psi(\zeta_{\bm \ell})}{Q(\zeta_{\bm \ell})^{\nu-1}}  +\sum_{\bm k\in\Lambda} q_{\bm k} c_{\bm k}.} Therefore, the dual problem of (\ref{opt3}) is
 \begin{align}\label{dual1}
& \hat Q=\underset{Q\in\Bc_+(N)}{\mathrm{argmin}}  \; J_\nu(Q) 
\end{align}
where $\Bc_+(N)=\{Q\in\Bc(N) \hbox{ s.t. } Q(\zeta_{\bm \ell})>0\, \forall\, \bm \ell \in\mathbb Z^2_N\}$ and $\Bc(N)$ is the set of all multidimensional trigonometric polynomials in $Q(\zeta_{\bm \ell})$ associated with the index set $\Lambda$.
\teo Assume that $\Lambda$ is such that $2n_j<N_j$ for $j=1,2$. Let $\{c_{\bm k}, \; \bm k \in\Lambda\}$ be the moments corresponding to $\Phi_o(\zeta_{\bm \ell})>0$, i.e. the original image. Moreover, we assume that $\Psi(\zeta_{\bm \ell})>0$ is fixed. Then, problem (\ref{dual1}) does admit a unique solution.\eteo
\proof It is not difficult to see that the second variation of $J_\nu(Q)$ along $\delta Q(\zeta_{\bm \ell})\in\Bc(N)$ is 
\al{\delta^2J_\nu(Q,\delta Q)=\frac{\nu}{|N|}\sum_{\bm \ell\in\mathbb Z_N^2} \frac{\Psi(\zeta_{\bm \ell})}{Q(\zeta_{\bm \ell})^{\nu+1}}\delta Q(\zeta_{\bm \ell})^2\geq 0.} If $\delta^2J(Q,\delta Q)=0$, then we have $\delta Q(\zeta_{\bm \ell})=0$ because $\Psi(\zeta_{\bm \ell})>0$ and $Q(\zeta_{\bm \ell})>0$. Since $2n_j<N_j$ with $j=1,2$, we have that   $\delta Q(\zeta_{\bm \ell})=0$ implies $\delta Q(\zeta_{\bm \ell})$ is the zero polynomial \cite[Lemma1]{7403052}. We conclude that $J_\nu$ is strictly convex in $\Bc_+(N)$. Thus, if $J_\nu$ admits minimum, then the latter is also unique. 

Next, we show that we can restrict the search of the minimum over a compact set in $\Bc_+(N)$. Since $J_\nu$ is continuous over $\Bc_+(N)$, by the Weierstrass theorem, $J_\nu$ does admit minimum which is also unique. We proceed to prove such a restriction is possible. Let $Q_j(\zeta_{\bm \ell})\in\Bc_+(N)$, with $j\in\mathbb N$, such that $Q_j(\zeta_{\bm \tilde \ell})\rightarrow \infty$ for some $\tilde{\bm \ell} \in\mathbb Z_N^2$ as  $j\rightarrow\infty$. Notice that 
\al{J_\nu(Q_j)&\geq \sum_{\bm k\in\Lambda} q_{\bm k} c_{\bm k}=\sum_{\bm k\in\Lambda} q_{\bm k} \frac{1}{|N|}\sum_{\bm \ell\in\mathbb Z_N^2} \zeta_{\bm \ell}^{\bm k} \Phi_o(\zeta_{\bm \ell})\nn\\
& =\frac{1}{|N|}\sum_{\bm \ell\in\mathbb Z_N^2}\sum_{\bm k\in\Lambda} q_{\bm k}  \zeta_{\bm \ell}^{\bm k} \Phi_o(\zeta_{\bm \ell})\nn\\&=\frac{1}{|N|}\sum_{\bm \ell\in\mathbb Z_N^2}Q(\zeta_{\bm \ell}) \Phi_o(\zeta_{\bm \ell})\nn \\ &\geq \frac{\min_{\bm \ell\in\mathbb Z_N^2} \Phi_o(\zeta_{\bm \ell}) }{|N| }
\sum_{\bm \ell\in\mathbb Z_N^2}Q(\zeta_{\bm \ell})\nn\\ &\geq  \frac{\min_{\bm \ell\in\mathbb Z_N^2} \Phi_o(\zeta_{\bm \ell}) }{|N|}
Q(\zeta_{\tilde{\bm \ell}}) \rightarrow \infty\nn} where we exploited the facts that $\Phi_o(\zeta_{\bm \ell})>0$ and \al{c_{\bm k}=\frac{1}{|N|}\sum_{\bm \ell\in\mathbb Z_N^2} \zeta_{\bm \ell}^{\bm k} \Phi_o(\zeta_{\bm \ell}), \; \; \bm k\in\Lambda.\nn}
Thus, we can restrict the search of $Q$ over the set $\Bc_1(N)=\{Q\in\Bc(N) \hbox{ s.t. } 0< Q(\zeta_{\bm \ell})\leq \mu \; \forall \,\bm\ell\in\mathbb Z_N^2\}$ for some $	\mu>0$. Next, let $Q_j(\zeta_{\bm \ell})\in\Bc_1(N)$, with $j\in\mathbb N$, be a sequence converging to $\bar Q(\zeta_{\bm \ell})$ for which $\bar Q(\zeta_{\tilde{\bm \ell}})=0$ for some $\tilde{\bm \ell}\in\mathbb Z^2_N$. Then $J_\nu(Q_j)\rightarrow \infty$ as $j\rightarrow \infty$. Hence, we can restrict the search of $Q$ over the bounded and closed (and thus compact) set $\Bc_2(N)=\{Q\in\Bc(N) \hbox{ s.t. } \gamma \leq Q(\zeta_{\bm \ell})\leq \mu \; \forall \,\bm\ell\in\mathbb Z_N^2\}$ for some $\gamma>0$.  
\qed
\rem In the primal problem (\ref{opt2}) the optimization variable $\Phi$ is defined over a grid of dimension $(2N_1+1)\times (2N_2+1)$, while in the dual problem (\ref{dual1}) the optimization variable $Q$ is defined over a grid of dimension $(2n_1+1)\times (2n_2+1)$. Since $n_j\ll N_j$, then it is more convenient to solve (\ref{dual1}) rather than (\ref{opt2}).\erem

The above result can be extended also for the case $\nu=1$ and $\nu=\infty$ as follows.  
\teo Consider the problem (\ref{opt2}) with $\nu=1$ and the same hypotheses of before. Then, the optimal solution is unique and such that $\Phi(\zeta_{\bm \ell})=\Psi(\zeta_{\bm \ell})/Q(\zeta_{\bm \ell})$ with $Q\in\Bc_+(N)$. The latter is given by the dual problem which is equivalent to minimize \al{J_1(Q)=\frac{1}{|N|}\sum_{\bm \ell\in\mathbb Z_N^2}\Psi(\zeta_{\bm \ell})\log\left(\frac{\Psi(\zeta_{\bm \ell})}{Q(\zeta_{\bm \ell})}\right)+\sum_{\bm k\in\Lambda} q_{\bm k} c_{\bm k}. } The latter does admit a unique solution.\eteo

\teo Consider the problem (\ref{opt2}) with $\nu=\infty$ and the same hypotheses of before. Then, the optimal solution is unique and such that $\Phi(\zeta_{\bm \ell})=\Psi(\zeta_{\bm \ell})\exp\left(-Q(\zeta_{\bm \ell}\right))$ with $Q\in\Bc(N)$. The latter is given by the dual problem which is equivalent to minimize \al{J_\infty(Q)=\frac{1}{|N|}\sum_{\bm \ell\in\mathbb Z_N^2}\Psi(\zeta_{\bm \ell})\exp\left(-Q(\zeta_{\bm \ell})\right)+\sum_{\bm k\in\Lambda} q_{\bm k} c_{\bm k}. } The latter does admit a unique solution.\eteo

\begin{figure*}
	\centering
	\begin{subfigure}{.32\textwidth}
		\centering
		\includegraphics[scale=0.3]{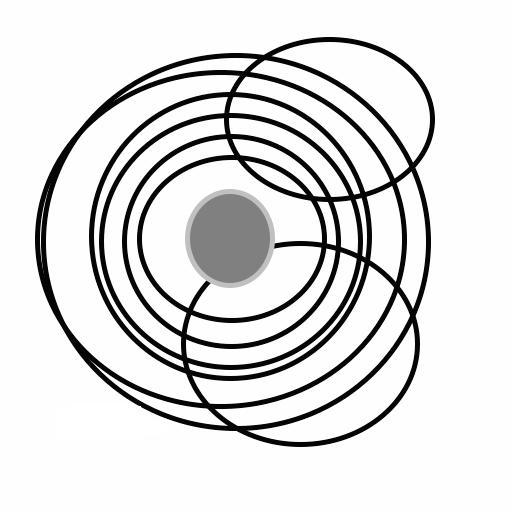} 
		
	\end{subfigure}
	\begin{subfigure}{.32\textwidth}
		\centering
		\includegraphics[scale=0.3]{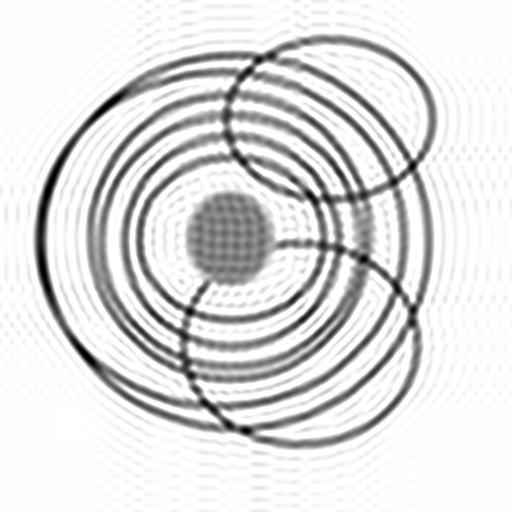}  
		
	\end{subfigure}
	\begin{subfigure}{.32\textwidth}
		\centering
		\includegraphics[scale=0.3]{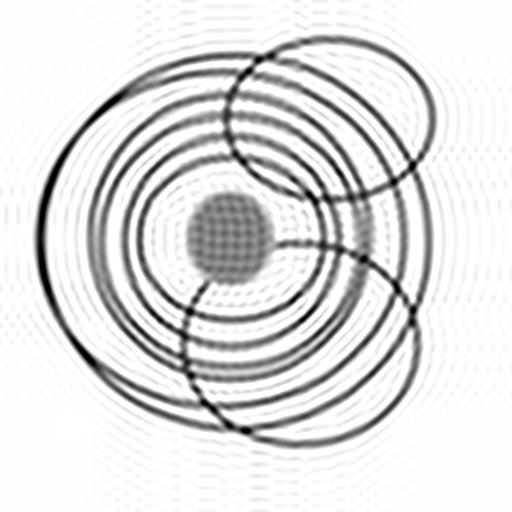}  
		
	\end{subfigure}
	\begin{subfigure}{.32\textwidth}
		\centering
		\includegraphics[scale=1.1,trim=50 10 0 17,clip]{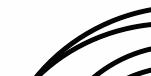} 
		\caption{Original image.} 
		\label{fig:original_circle}
	\end{subfigure}
	\begin{subfigure}{.32\textwidth}
		\centering
	\includegraphics[scale=1.1,trim=50 10 0 17,clip]{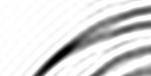}  
	\caption{$\nu=1$, $PSNR=15.529 \text{ dB}$.}
		\label{fig:circlesv=1}
	\end{subfigure}
	\begin{subfigure}{.32\textwidth}
		\centering
	\includegraphics[scale=1.1,trim=50 10 0 17,clip]{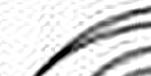}  
	\caption{$\nu=\infty$, $PSNR = 16.202 \text{ dB}$.}
		\label{fig:circlesv=inf}
	\end{subfigure}
	\caption{Effect of the parameter $\nu$ using the Alpha-divergence with $cr\approx0.97$.}
\end{figure*}

\rem In the case that $\nu=1$ with $\Psi(\zeta_{\bm \ell})=1$, i.e. the uniform prior, the objective function in problem (\ref{opt2}) is equivalent to the entropy functional in (\ref{def_entropy}).\erem Since the optimal solution to (\ref{opt2}) depends on $\nu$, we can compress the original image also including the value of $\nu$ for which we have the best image reconstruction from the moments. The performance of the image reconstruction is measured by the peak signal-to-noise ratio (PSNR) which is the typical index adopted in image processing, \cite{psnr}:
\al{ PSNR= 10\log_{10}\left(\frac{MAX_{X^o}^2}{MSE}\right).}
where $MSE = \sum_{i,j}(X^o_{i,j}-\hat{X}_{i,j})^2/{p_1p_2}$,
$X^o$ is the original image $\Phi_o$ represented in matrix form, $\hat{X}$ is the reconstructed image $\hat \Phi$ represented in matrix form, and $MAX_{X^o}$ is the maximum possible pixel value of the original image, which in this case is equal to 1 since $X^o\in[0,1]^{p_1\times p_2}$.
 Algorithm \ref{algo2} and \ref{algo3} show the compression and the reconstruction procedure, respectively. \begin{algorithm}
 \KwInput{$\Phi_o$, $\Lambda$ (i.e. compression rate)}
 \KwOutput{$\{c_{\bm k}, \; \bm k\in \Lambda\}$, $\nu^{opt}$}
 Compute $c_{\bm k}$, $\bm k\in\Lambda$ as in (\ref{moment_cond})\;  
 Let $\nu_1\ldots \nu_m$ be a set of candidates for $\nu$\;
 \For{k=1:m}{
Compute $\hat \Phi_{\nu_k}$ solution to (\ref{opt2}) with $\nu=\nu_k$\;
Let $PSNR_k$ be the PSNR between $\Phi_o$ and $\hat \Phi_{\nu_k}$}
Let $k^{opt}$ be the one minimizing $PSNR_k$\;
Set $\nu^{opt}=\nu_{k^{opt}}$
  \caption{Image compression}\label{algo2}
\end{algorithm}
\begin{algorithm}
 \KwInput{$\{c_{\bm k}, \; k\in \Lambda\}$, $\nu^{opt}$}
 \KwOutput{$\hat\Phi$ (i.e. the reconstructed image)}
Compute $\hat \Phi$ solution to (\ref{opt2}) with $\nu=\nu^{opt}$\;
  \caption{Image reconstruction}\label{algo3}
\end{algorithm} We apply our procedures to the original image in Figure \ref{fig:original_circle} ($p_1=p_2=512$). 

\begin{figure}\centering
\includegraphics[width=0.8\columnwidth]{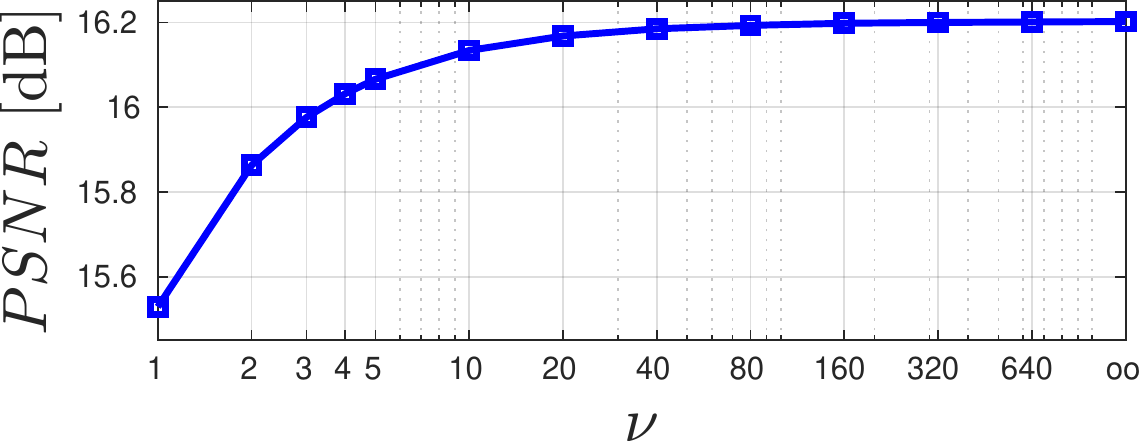}
\caption{ Considered values for $\nu$ and the corresponding PSNR between $\Phi_o$ and $\hat \Phi_\nu$.}
\label{plotnu}
\end{figure}

For simplicity we choose $\Psi(\zeta_{\bm \ell})=1$, i.e. the uniform prior. Here, $n_1=n_2=85$ meaning that the compression rate is $cr=1-((n_1+1)(n_2+1)+1)/p_1p_2\approx 0.97$. The set of the candidates for $\nu$ is depicted in Figure \ref{plotnu} with the corresponding PSNR. 
Clearly, $\nu^{opt}=\infty$. Figures \ref{fig:circlesv=1} and \ref{fig:circlesv=inf} show the reconstructed images for $\nu=1$ and $\nu=\infty$, respectively. Although the images seem to be similar, the image loses in focus for $\nu=1$: see the zoom in Figures \ref{fig:original_circle}, \ref{fig:circlesv=1} and \ref{fig:circlesv=inf}. We conclude that  it is possible to highlight specific details of the image, depending on the context of application, by properly selecting the objective function.

\section{Design of the prior}\label{sect_prior}
Recall that the prior $\Psi$ embeds some information about the image. In what follows, we show how the prior can be designed in two different scenarios.  

\subsection{Hybrid compression}
The image compression method proposed in \cite{7403052} designs $\Psi$ in such a way the cepstral moments in $\Lambda\setminus \{0\}$ of the original image are approximately matched. In this scenario both the moments in (\ref{moment_cond}) and the cepstral moments constitute the compressed image. Roughly speaking, $\Psi$ is extracted from the cepstral moments. In what follows we propose to extract $\Psi$ form the singular values decomposition (SVD) of the original image. More precisely, let $X\in\mathbb R^{p_1\times p_2} $ be the nonnegative matrix representing the original image. Let $X=UDV^T$ be the SVD of $X$ where $UU^T=I_{p_1}$, $VV^T=I_{p_2}$ and $D\in\mathbb{R}^{p_1\times p_2}$ is the matrix containing the singular values $d_k$, $k=1\ldots \min\{p_1,p_2\}$. We assume that $d_k\geq d_l$ with $k\leq l$. Then, we choose $r\ll\min\{p_1,p_2\} $ and we consider the following low-rank approximation $X_{r}=U \bar D_rV^T$
where $\bar D_r$ is obtained from $D$ by substituting $d_k$ with $0$ for $k>r$. Then, we define the prior as $\Psi (\zeta_{\bm \ell}) =\exp(Y_{\ell_1+1,\ell_2+1})$ where $Y\in\mathbb R^{N_1\times N_2}$ is the symmetric mirroring of $X_r$.
Therefore, to construct $\Psi$ we need to store the matrices 

 {\small \al{ & M_{1,r}=U \left[\begin{array}{c}D_r^{1/2}  \\ 0 \end{array}\right]\in\Rs^{p1\times r}, \, M_{2,r}^T=V\left[\begin{array}{cc}D_r^{T/2}  \\ 0^T \end{array}\right]\in\Rs^{p_2\times r} \nn}}
 
\noindent where $D_r^{1/2}=\mathrm{diag}(\sqrt{d_1}\ldots \sqrt{d_{r}})$, i.e. we have $(p_1+p_2)r$ parameters. Accordingly, the compressed image is given by the moments set $\{ c_{\bm k}, \; \bm k\in\Lambda\}$, $M_{1,r}$ and $M_{2,r}$.  Accordingly, the compression is hybrid in the sense we combine two difference strategies for this aim: moments and SVD. The compression rate is \al{\label{cr_form}cr=1-((p_1+p_2)r+(n_1+1)(n_2+1))/p_1 p_2.} Then, to reconstruct  the image we extract the prior, say $\Psi_r$, from $M_{1,r}$ and $M_{2,r}$ as above then we solve problem (\ref{opt2}) with $\Psi=\Psi_r$ where for simplicity $\nu$ is kept fixed. Clearly,  the optimal solution to problem (\ref{opt2}), say $\hat \Phi_r$, depends on the choice of $r$. Therefore, once $cr$ is fixed, we can select $r$, $n_1$ and $n_2$, compatible with the chosen $cr$ through (\ref{cr_form}), leading the to best image extraction in terms of PSNR. For simplicity, assume that $n_1=n_2$, then given $cr$ and $r$ we have that 
 \al{\label{eqn}n_1=n_2=\mathrm{round}(\sqrt{(1- cr)p_1p_2-(p_1+p_2)r}).} Algorithm \ref{algo1} shows the compression procedure. \begin{algorithm}
 \KwInput{$\Phi_o$, $\Lambda$ (i.e. compression rate)}
 \KwOutput{$\{c_{\bm k}, \; \bm k\in \Lambda\}$, $M_1$, $M_2$}
Compute $c_{\bm k}$, $\bm k\in\Lambda$ as in (\ref{moment_cond})\;  
 \For{$r=0:r_{max}$}{
Compute $M_{1,r}$, $M_{2,r}$ from the SVD\;
Extract the prior $\Psi_r$ from $X_r=M_{1,r}M_{2,r}$\;
Set $n_1$, $n_2$ by (\ref{eqn})\;
Compute $\hat \Phi_r$ solution to (\ref{opt2}) with $\Psi=\Psi_r$\;
Let $PSNR_r$ be the PSNR between $\Phi_o$ and $\hat \Phi_{r}$}
Let $r^{opt}$ be the one minimizing $PSNR_r$\;
Set $M_1=M_{1,r^{opt}}$, $M_2=M_{2,r^{opt}}$\;
  \caption{Image compression}\label{algo1}
\end{algorithm} Notice that $r_{max}$ is the largest integer for which (\ref{eqn}) leads to a real solution for $n_1$ and $n_2$. Finally, the extraction procedure is similar to Algorithm \ref{algo2}: the unique differences are that: $\nu$ is not a parameter of the compressed image; the prior now is extracted from $M_1$ and $M_2$. 

\begin{figure}\centering
\includegraphics[width=0.8\columnwidth]{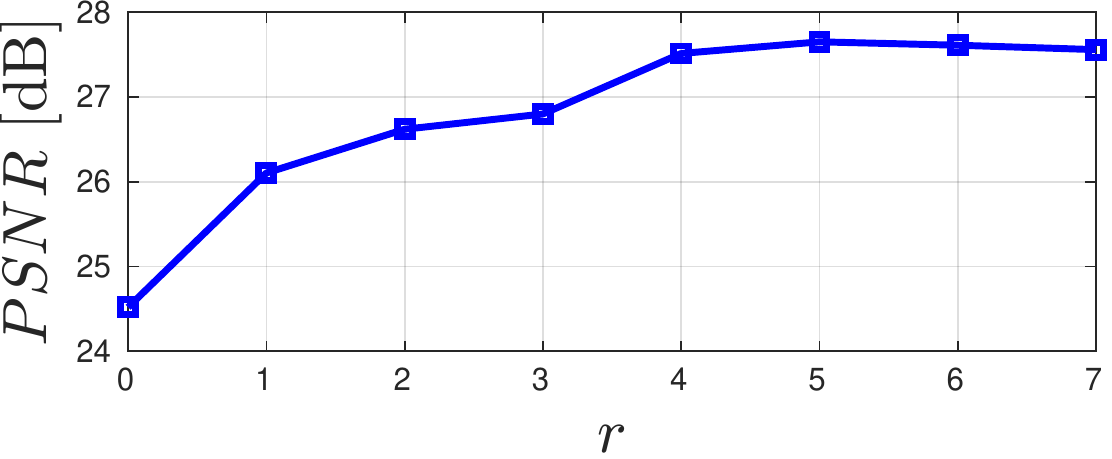}
\caption{ Considered values for $r$ and the corresponding PSNR between $\Phi_o$ and $\hat \Phi_r$.}\label{fig_r}
\end{figure}

We tested our procedure with the image taken from \href{https://www.archdaily.com/932368/the-platform-office-building-space-encounters}{\sl https://www.archdaily.com/932368/the-platform-office-building-space-encounters}, with $p_1=p_2=512$. Choosing $cr\approx0.97$, we have $r_{max}=7$. Figure \ref{fig_r} shows the  PSNR of all the possible values of $r$, for $\nu=\infty$ fixed; clearly we have $r^{opt}=5$. By applying the compression algorithm with a high $r$ number, the resulting image provides better details of the building structure. By contrast, in order to obtain a well detailed version of the machinery in the bottom left-hand corner, a large number of moments $n_1,n_2$ is preferable (thus, a small value of $r$). Clearly, the optimal balance among $n_1,n_2$ and $r$ strongly depends on the type of the considered image. In particular, the more the details are disposed according to a vertical and horizontal grid pattern, the more the SVD is efficiently capable of representing them and therefore a large value of $r$ will be likely preferable. In the case of the considered image, many details are placed along an orthogonal-grid pattern and therefore the result corresponding to a small value of $r$ is satisfying. Nevertheless, since the moments are in general much less expensive in terms of stored memory, it suffices to slightly increase the number of stored moments $n_1,n_2$ (at the price of minimally decreasing the number of stored singular values $r$) to obtain a reasonable compromise for the overall image.

Finally, Algorithm \ref{algo1} can be easily extended to the case in which $\nu$ is optimized as in Section \ref{sect_obj}. However, such an extension increases the computational burden because we have to construct a two-dimensional grid, indexed by $\nu$ and $r$, in order to select the best reconstructed image.

\begin{figure*}
\centering
\begin{subfigure}{.32\textwidth}
\centering
\includegraphics[scale=0.5]{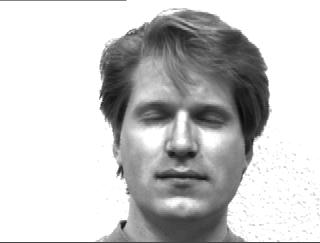} 
\end{subfigure}
\begin{subfigure}{.32\textwidth}
\centering
\includegraphics[scale=0.5]{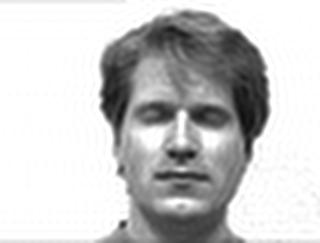}  
\end{subfigure}
\begin{subfigure}{.32\textwidth}
\centering
\includegraphics[scale=0.5]{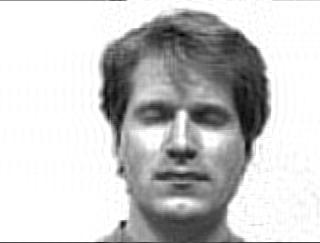}  
\end{subfigure}
\begin{subfigure}{.32\textwidth}
\centering
\includegraphics[scale=1.3,trim=0 0 0 3,clip]{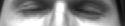} 
\caption{Original image.} \label{fig:face_original}
\end{subfigure}
\begin{subfigure}{.32\textwidth}
\centering
\includegraphics[scale=1.3,trim=0 0 0 3,clip]{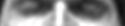}  
\caption{Uniform prior, $PSNR=26.835\text{ dB}$.}\label{fig:face_r0}
\end{subfigure}
\begin{subfigure}{.32\textwidth}
\centering
\includegraphics[scale=1.3,trim=0 0 0 3,clip]{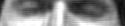}  
\caption{Low-rank prior, $PSNR=29.506\text{ dB}$.}\label{fig:face_r15}
\end{subfigure}
\caption{Image extraction using a uniform prior versus a similar image as prior. The entire image is reported above and a zoom on the eyes below.}
\label{fig:faces}
\end{figure*}
\subsection{Compression of a database of similar images}
\text{\the\columnwidth} Assume that we have a database of similar images. We compress these images by computing the moments set $\{c_{\bm k}, \; \bm k\in \Lambda\}$ for each image. In such a scenario, we can extract these images by using a common prior $\Psi$. The latter can be obtained, for instance, from a rough approximation of one of these images, e.g. a low rank approximation from the SVD of one image of the database.

\begin{figure}
\centering
\includegraphics[scale=0.45]{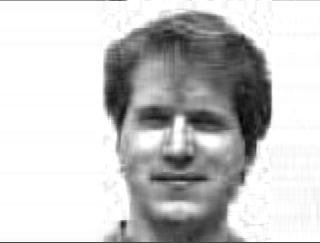} 
\caption{Low-rank approximation of the similar image used as prior. The approximation is obtained by means of SVD with rank equal to $15$.} 
\label{fig:face_similar}
\end{figure}

We apply this approach to the {\em yalefaces} dataset, \cite{YALE}. 
More precisely, we consider the original image in Figure \ref{fig:face_original}, here $p_1=243$ and $p_2=320$. We compress the image by computing $n_1=53$ and $n_2=70$ moments. Then, we extract the image by solving problem (\ref{opt2}) with $\nu=1$ and the prior is chosen in two different ways. In the first one we use $\Psi(\zeta_{\bm \ell})=1$, i.e. the uniform prior.  In the second one, we choose as prior the low rank approximation of another (similar) image in the dataset, see Figure \ref{fig:face_similar}. 

Figures \ref{fig:face_r0} and \ref{fig:face_r15} show the extracted images. The use of the second prior increases the performances of the extraction algorithm both in terms of PSNR and in terms of details. In fact it can be seen that the extracted image in Figure \ref{fig:face_r15} is much sharper than the one using the uniform prior. Indeed, this is more evident by looking at the the details around the eyes, see Figures \ref{fig:face_original}, \ref{fig:face_r0} and \ref{fig:face_r15}. Finally, we have obtained similar results by changing the parameter $\nu$.

\section{Conclusions}\label{sect_concl}
We have revisited the image compression problem based on multidimensional circulant covariance extension problem which has been proposed in \cite{7403052}. In particular, we have explored the possibility to change the objective function and the prior in order to refine the image reconstruction. These are preliminary ideas that could be further extended. For instance, the compression method for a database of similar images could be adapted to efficiently compress video sequences. Indeed, in such a scenario, adjacent frames are reasonably similar one another. Accordingly, the extracted frame at time $t$ could be used as prior to extract the frame at time $t+1$.   


\end{document}